\documentclass[9pt]{article}
\usepackage{amsmath}
\usepackage{amsfonts}
\usepackage{amssymb}

\newtheorem{lemma}{Lemma}[section]

\newtheorem{theorem}{Theorem}[section]
\newtheorem{corollary}{Corollary}[section]
\newtheorem{remark}{Remark}[section]
\newtheorem{definition}{Definition}[section]
\newtheorem{example}{Example}[section]

\let\Section=\section
\def\section{\setcounter{equation}{0}\Section}
\sf

\begin{document}
\title{ Bound states for a stationary nonlinear Schr\"{o}dinger-Poisson system with
sign-changing potential in $\mathbb{R}^3$
\thanks{This work was supported by
NSFC and CAS-KJCX3-SYW-S03. $^\dag$ Corresponding
author.\hfil\break\indent {\it 2000 Mathematical Subject
Classification}: 35J60. \hfil\break\indent{\it Key words}: Schr\"
odinger-Poisson system, sub-supersolutions, supercritical Sobolev
exponent, sign-changing potential, bound state.
 \hfil\break \indent To appear in Acta Math. Sci., 29B(2009), No.4, p. 1095- }}
 \author{Yongsheng Jiang and Huan-Song Zhou$^\dag$}
 \date{}
\maketitle

\noindent {\sc Abstract}: We study the following
Schr\"{o}dinger-Poisson system
\begin{equation}\nonumber
(P_\lambda)\left\{\begin{array}{ll}
 -\Delta u + V(x)u+\lambda \phi (x) u =Q(x)u^{p}, x\in \mathbb{R}^3 \\
 -\Delta\phi = u^2, \lim\limits_{|x|\rightarrow +\infty}\phi(x)=0,
 u>0,
 \end{array}\right.
\end{equation}
where  $\lambda\geqslant0$ is a parameter, $1 < p < +\infty$,
$V(x)$ and $Q(x)$ are sign-changing or non-positive functions in $
L^{\infty}(\mathbb{R}^3)$. When $V(x)\equiv Q(x)\equiv1$, D.Ruiz
\cite{RuizD-JFA} proved that ($P_\lambda$) with $p\in(2,5)$ has
always a positive radial solution, but ($P_\lambda$) with
$p\in(1,2]$ has solution only if $\lambda>0$ small enough and no
any nontrivial solution if $\lambda\geqslant\frac{1}{4}$. By using
sub-supersolution method, we prove that there exists $\lambda_0>0$
such that ($P_\lambda$) with $p\in(1,+\infty)$ has always a bound
state ( $H^1(\mathbb{R}^3)$ solution) for
$\lambda\in[0,\lambda_0)$ and
 certain functions $V(x)$ and $Q(x)$ in $
L^{\infty}(\mathbb{R}^3)$. Moreover, for every
$\lambda\in[0,\lambda_0)$, the solutions $u_\lambda$ of $\rm
(P_\lambda)$ converges, along a subsequence, to a solution of
($P_0$) in $H^1$ as $\lambda \rightarrow 0$.

\section{Introduction}

  In this paper, we are concerned with the existence of positive solutions of the following nonlinear elliptic system
\begin{equation}\label{eq:1.1}
\left\{\begin{array}{ll}
 -\Delta u + V(x)u+\lambda\phi (x) u =Q(x)|u|^{p-1}u,\,\,\, x \in \mathbb{R}^3, \\
 -\Delta\phi = u^2,\,\,\lim\limits_{|x|\rightarrow +\infty}\phi(x)=0,
 \end{array}\right.
\end{equation}
where $\lambda>0$ is a parameter, $p\in(1,+\infty)$, $V(x)$ and
$Q(x)$ are functions in $L^{\infty}(\mathbb{R}^3)$. This kind of
problem is related to looking for solitary wave type solution of
nonlinear Schr\"{o}dinger equation for a particle in a
electromagnetic field \cite{d'AveniaP-ANonSt}, for more physical
background about this system we refer the reader to
\cite{VBenciDFortucto-TopMNA,KhalidBenmlih-EJDE,d'AveniaP-ANonSt,NIERF-PRSESM,RuizD-preprint,OScarSanchezJuanSoer-JSP,
WangZhengpingZhouHuan-Song-DisCDS} and the references therein.
Under variant assumptions on $V(x)$ and $Q(x)$, problem
(\ref{eq:1.1}) has been studied widely. For $Q(x)\equiv0$ and
$V(x)\equiv constant$, this problem was studied as an eigenvalue
problem in \cite{VBenciDFortucto-TopMNA} on bounded domain and in
\cite{KhalidBenmlih-EJDE,NIERF-PRSESM} on $\mathbb{R}^3$. For
$Q(x)\equiv1$ with $p\in(1,5)$, there has been quite a lot of
interest on problem (\ref{eq:1.1}) in recent years. For examples,
the existence of solutions to problem (\ref{eq:1.1}) with
$V(x)\equiv constant$ and $\lambda=1$ was obtained in
\cite{d'AveniaP-ANonSt} if $p\in(3,5)$ and in
\cite{D'AprileTMugnaiD-PRSESM} if $p\in[3,5)$, then in
\cite{RuizD-JFA} for $p\in(1,5)$ and $\lambda$ may not be equal to
1. Moreover, the existence of multiple solutions of (\ref{eq:1.1})
with $V(x)\equiv Q(x)\equiv1$ and $p\in(1,5)$ was proved by
Ambrosetti-Ruiz in
\cite{AmbrosettiAntonioRuizDavid-CCM,AmbrosettiA-MJM}. If $V(x)$
is not a constant, some existence results on problem
(\ref{eq:1.1}) were given in \cite{AzzolliniAandPomponioA-Jmaa}
for $Q(x)\equiv1$ with $p\in(3,5)$, then in
\cite{ZhaoLeigaZhaoFukun-Jmaa} for $p\in(2,3]$, and in
\cite{DSY,WangZhengpingZhouHuan-Song-DisCDS} for a general
nonlinear term $f(x,u)$. If $V$ and $Q$ are radial, positive, and
vanishing at infinity, the existence and nonexistence of solutions
to (\ref{eq:1.1}) were studied in \cite{MercuriC} for some $p$ in
$(1,5)$.  The results obtained in all the papers mentioned above
are based on variational methods, this leads to the restriction on
$p\in(1,5]$. For $V(x)\equiv Q(x)\equiv1$ and $\lambda=1$, it was
proved in \cite{D'AprileT-ANonSt,RuizD-JFA} that problem
(\ref{eq:1.1}) does not possess any nontrivial solution if
$p\leqslant2$ or $p\geqslant5$. What would happen if $V(x)$ and
$Q(x)$ are not equal to 1? Is it possible to get a solution of
problem (1.1) for all $p\in(1,+\infty)$?  To the authors'
knowledge, there seems no any results in this direction. In this
paper, we prove that for any $p\in(1,+\infty)$ and for certain
$V(x)$, there always exists $Q(x)$ such that problem
(\ref{eq:1.1}) has a positive solution if $\lambda>0$ small. As it
is known, if $p\in(1,+\infty)$, the variational approach is no
longer applicable and here we use sub-supersolution method
instead. But problem (\ref{eq:1.1}) is a coupled system, it seems
not easy to construct a reasonable sub- and supersolutions to
ensure the existence of a solution to the problem. Motivated by
the paper of Edelson-Stuart \cite{EdelsonALStuartCA-JDE} and based
on an estimate for the fundamental solution $\phi$ of the second
equation in (\ref{eq:1.1}), we get the desired sub- and
supersolutions of (\ref{eq:1.1}) for some kinds of $V(x)$, $Q(x)$
and $\lambda\geqslant0$ small. Therefore, by an iterative
procedure, we obtain a solution $u_\lambda$ of (\ref{eq:1.1}) for
each $\lambda\geqslant0$ small enough. In particular, our results
imply the existence of positive solution to the following single
equation
\begin{equation}\label{eq:1.2}
 -\Delta u + V(x)u =Q(x)|u|^{p-1}u,\,\,\, x \in \mathbb{R}^3,
\end{equation}
where $V(x)$ and $Q(x)$ are functions in $L^\infty(\mathbb{R}^3)$.
Moreover, we prove that, along a subsequence, the solutions
$u_\lambda$ of (\ref{eq:1.1}) for
$\lambda\in(0,-2(2\alpha-1)\Lambda)$ converges in
$H^1(\mathbb{R}^3)$ to a solution of (\ref{eq:1.2}), where
$\alpha>\frac{3}{4}$ and $\Lambda<0$ are given by $(H_3)$ and
$(H_4)$ below, respectively. (\ref{eq:1.2}) is essentially the
special case of (\ref{eq:1.1}) as $\lambda=0$, and it has been
studied by many authors, such as \cite{CostaDGTehraniH-CalVarPDE,
ChabrowskiJCostaDG-CPDE, ChenJQLiSJ-ManusMath}, etc. However, in
those papers, $p\in(1,5)$, $V(x)$ is assumed to be of the form
$\lambda h(x)$ and $Q(x)$ is required to have a negative limit as
$|x|\rightarrow+\infty$, $\int Q(x)\phi_1^{p+1}(x)dx<0$ or
$Q(x)\equiv1$, where $\phi_1>0$ is the first eigenfunction of the
problem (see e.g. \cite[Corollary
2]{AllegrettoWHuangYX-FunkcElvac})
\begin{equation}\label{eq:1.3}
 -\Delta u=\mu V(x) u, \ x\in \mathbb{R}^3 \mbox{ and } u \in D^{1,2}(\mathbb{R}^3).
\end{equation}
In some sense, our result on (\ref{eq:1.2}) also generalizes that
of
\cite{ChabrowskiJCostaDG-CPDE,ChenJQLiSJ-ManusMath,CostaDGTehraniH-CalVarPDE}.
Specially, in our case, $p\in(1,+\infty)$ is allowed, and we do
not require that
$\overline{\{x\in\mathbb{R}^3:Q(x)>0\}}\cap\overline{\{x\in\mathbb{R}^3:Q(x)<0\}}=\emptyset$.
See our Examples \ref{ex:1.1} and \ref{ex:1.2}.\\

Now, we give our assumptions on $V(x)$ and $Q(x)$.
\begin{description}
 \item $\rm (H_1)$  $V(x)$, $Q(x)$ are nonpositive, or
 sign-changing,
 functions in $C^{0,\gamma}_{loc}(\mathbb{R}^3,\mathbb{R})\cap L^{\infty}(\mathbb{R}^3
)$ with $\gamma\in(0,1)$.
 \item $\rm (H_2)$ There exists a
constant $V_\infty \geq 0$ such that
$\liminf\limits_{|x|\rightarrow\infty}V(x)=V_\infty$.
 \item $\rm
(H_3)$ 
$\Lambda:=\inf\{\int_{\mathbb{R}^3}|\nabla u|^2+V(x)u^2dx:u\in
H^1(\mathbb{R}^3),\int_{\mathbb{R}^3} u^2dx=1\} <0.$
\item $\rm (H_4)$ There exists $\alpha> {3}/{4}$ such that
\begin{equation}\label{Q}
 Q(x)\leqslant(1+r^2)^{\alpha(p-1)}[V(x)-H(r)]\ \text{for }
x\in\mathbb{R}^3,
\end{equation}
where $H(r)=2\alpha[(2\alpha-1)r^2-3](r^2+1)^{-2}$ and $r=|x|$.
Moreover,
 \begin{equation} \label{V}
  \lim\limits_{|x|\rightarrow +\infty} |x|^2V(x) > 2\alpha
(2\alpha-1) \mbox{ if } V_\infty =0 \mbox{ and } p\geq 1+1/\alpha.
\end{equation}
\end{description}
\begin{remark}\label{re:1.1}
Note that (\ref{V}) is only used to ensure that $Q(x)$ with
property (\ref{Q}) is not $-\infty$. Condition $\rm (H_4)$ can be
slightly weakened by assuming that
\begin{description}
\item $\rm (H_4)'$ There exist $\alpha>\frac{3}{4}$, $\theta>0$ and
$a>0$ such that
\begin{equation}\nonumber
Q(x)\leqslant
a^{1-p}(1+\theta^2r^2)^{\alpha(p-1)}[V(x)-\theta^2H(\theta r)],\  \
x\in\mathbb{R}^3.
\end{equation}
\end{description}
\end{remark}
\begin{remark}\label{re:1.2}
For $H(r)$ given by $\rm(H_4)$, let
$r_0=\sqrt{\frac{3}{2\alpha-1}}$, we see that $H(r_0)=0$, $H(r)<0$
for $r\in(0,r_0)$ and $H(r)>0$ for $r>r_0$.
\end{remark}
Here are two examples on our assumptions. Example \ref{ex:1.1}
satisfies $\rm(H_1)-(H_4)$ and the assumptions of
\cite{ChabrowskiJCostaDG-CPDE}. Example \ref{ex:1.2} satisfies
$\rm(H_1)-(H_4)$, but does not satisfy  the assumptions of
\cite{ChabrowskiJCostaDG-CPDE,CostaDGTehraniH-CalVarPDE}.
\begin{example}\label{ex:1.1}
Let $\alpha>\frac{3}{4}$,  $b>1$ and $\beta>0$. For $H(r)$ given
by $\rm(H_4)$, let $V(x)=bH(r)$ and
$Q(x)\leqslant(b-1)(\frac{2\alpha+2}{2\alpha-1})^{\alpha(p-1)}H(r)$
with $Q(x)\in C^{0,\gamma}_{loc}(\mathbb{R}^3,\mathbb{R})\cap
L^{\infty}(\mathbb{R}^3 )$ and
$\lim\limits_{|x|\rightarrow\infty}Q(x)=-\beta$. Then, $\rm(H_1)$
and $\rm(H_2)$ with $V_\infty=0$ are satisfied. For $r_0$ given by
Remark \ref{re:1.2}, by taking $\varphi(x) \geq 0$ and $\varphi
\in C^{\infty}_0(B_{r_0}(0))\setminus \{0\}$ and we see that
$\rm(H_3)$ is satisfied if $b>1$ large enough. By Remark
\ref{re:1.2} and a directly computation shows that  $\rm(H_4)$ is
also satisfied. Moreover, for any $p>1$, it follows from
 $\lim\limits_{|x|\rightarrow\infty}Q(x)=-\beta$ that there is $\beta_0>0$ such that
$\int_{\mathbb{R}^3}Q(x)\phi_1^{p+1}<0$ if $\beta>\beta_0$, where
$\phi_1>0$ is the first eigenfunction of (\ref{eq:1.3}). Hence,
the conditions of \cite{ChabrowskiJCostaDG-CPDE} are also
satisfied.
\end{example}
\begin{example}\label{ex:1.2}
In Example \ref{ex:1.2}, we take
$Q(x)=(b-1)(\frac{2\alpha+2}{2\alpha-1})^{\alpha(p-1)}H(r)$, and
now $\beta=0$. Then we still have that $\rm(H_1)-(H_4)$ are
satisfied for $b>1$ large. But  the condition on $Q(x)$ in
\cite{ChabrowskiJCostaDG-CPDE,CostaDGTehraniH-CalVarPDE} cannot be
satisfied because here we have that
$\int_{\mathbb{R}^3}Q(x)\phi_1^{p_0+1}>0$ for some $p_0>1$. In
fact, since $s>3$ we know that $V(x)\in
L^{\frac{s}{2}}(\mathbb{R}^3)\cap L^{\infty}(\mathbb{R}^3)$, it
follows from \cite[Corollary 2]{AllegrettoWHuangYX-FunkcElvac}
that the first eigenvalue $\mu_1$ of (\ref{eq:1.3}) is positive
and it has an positive eigenfunction $\phi_1\in
D^{1,2}(\mathbb{R}^3)$. Hence $\int
V(x)\phi_1^2=\frac{1}{\mu_1}\int|\nabla \phi_1|^2>0$ by
(\ref{eq:1.3}), this implies that
$\int_{\mathbb{R}^3}Q(x)\phi_1^{2}>0$. Moreover, by \cite[Theorem
8.17]{DGilbargNStrudinger} we see that $\phi_1 \in
L^\infty(\mathbb{R}^3)$. Then the dominated convergence theorem
shows that $\int_{\mathbb{R}^3}Q(x)\phi_1^{p+1} \rightarrow
\int_{\mathbb{R}^3}Q(x)\phi_1^{2}>0$ as $p\rightarrow 1^+$. So,
there is some $\delta>0$ such that
$\int_{\mathbb{R}^3}Q(x)\phi_1^{p+1}>0$ for $p\in (1, 1+\delta)$.
\end{example}

Finally, we give the main results of the paper.
\begin{theorem}\label{th:1.1}
For any $p\in(1,+\infty)$, suppose that $\rm(H_1)$ to $\rm(H_4)$
are satisfied and $\Lambda<V_\infty$. Then, problem (\ref{eq:1.1})
has at least a positive solution $u_\lambda\in
C_{loc}^{2,\gamma}(\mathbb{R}^3)\cap W^{2,q}(\mathbb{R}^3)$ for
any $\lambda\in[0,-2(2\alpha-1)\Lambda)$ and all
$q\in[2,+\infty)$. Moreover, $u_\lambda$ is a bound state with
\begin{equation}\label{eq:1.5}
0<u_\lambda\leqslant \frac{1}{(1+|x|^2)^\alpha},
\end{equation}
and
\begin{equation}\label{eq:1.6}
\|u_\lambda\|_{W^{2,q}(\mathbb{R}^3)}\leqslant C,
\end{equation}
where $C$ is a constant independent of $\lambda$.
\end{theorem}
In particular, if $\lambda=0$ in (\ref{eq:1.1}), Theorem
\ref{th:1.1} implies that
\begin{corollary}\label{co:1.1}
Under the assumptions of Theorem 1.1, problem (\ref{eq:1.2})
possesses a positive solution $u\in
C_{loc}^{2,\gamma}(\mathbb{R}^3)\cap W^{2,q}(\mathbb{R}^3)$ for all
$q\in[2,+\infty)$.
\end{corollary}
\begin{theorem}\label{th:1.2} For each $\lambda\in
(0,-2(2\alpha-1)\Lambda)$, let $u_{\lambda}$ denote a solution of
problem (\ref{eq:1.1}) obtained by Theorem \ref{th:1.1}, then
there exists $u_0\in C_{loc}^{2,\gamma}(\mathbb{R}^3)\cap
W^{2,q}(\mathbb{R}^3)$ for all $q\in[2,+\infty)$ such that, along
a subsequence,
\begin{equation}\nonumber
\|u_{\lambda}-u_0\|_{W^{1,2}(\mathbb{R}^3)}\rightarrow0,\ \text{as}\
\lambda\rightarrow0,
\end{equation}
and $u_0$ is a positive solution of (\ref{eq:1.2}).
\end{theorem}
\begin{remark} \label{rek1.3}
We believe that our methods for proving the above results work
also for the general case $N>3$.
\end{remark}

Throughout this paper, we denote the usual norm of
$L^q(\mathbb{R}^3)$ and $W^{2,q}(\mathbb{R}^3)$ for
$q\in[1,+\infty]$, respectively, by $|\cdot|_q$ and
$\|\cdot\|_{2,q}$.

\section{ Subsolution and Supersolution}
The aim of this section is to construct a subsolution and a
supersolution of problem (\ref{eq:1.1}). Based on these sub- and
supersolutions Theorem \ref{th:1.1} is proved in Section 3. We
begin this section by giving our definitions of sub- and
supersolutions for system (\ref{eq:1.1}).
\begin{definition}\label{de:2.1}
A positive function $\psi(x)\in C^2(\mathbb{R}^3)$ is said to be a
{\bf supersolution} of (\ref{eq:1.1}) if
\begin{equation}\label{eq:2.1}
-\Delta\psi(x)+V(x)\psi(x)+\lambda\phi(x)\psi(x)\geqslant
Q(x)\psi^{p}(x),\  x\in\mathbb{R}^3,
\end{equation}
with $\phi(x)$ satisfies
\begin{equation}\label{eq:2.2}
-\Delta\phi(x)= u^{2}(x),\
\lim\limits_{|x|\rightarrow+\infty}\phi(x)=0,
\end{equation}
$\ \text{for}\ u\in W^{2,2}(\mathbb{R}^3)\ \text{and}\
0<u(x)\leqslant\psi(x)\ \text{on}\ \mathbb{R}^3.$ A positive
function $\varphi(x)\in C^2(\mathbb{R}^3)$ is said to be a {\bf
subsolution} of (\ref{eq:1.1}) if the opposite inequality to
(\ref{eq:2.1}) is satisfied by $\varphi(x)$, that is
\begin{equation}\label{eq:2.3}
-\Delta\varphi(x)+V(x)\varphi(x)+\lambda\phi(x)\varphi(x)\leqslant
Q(x)\varphi^{p}(x),\  x\in \mathbb{R}^3,
\end{equation}
and (\ref{eq:2.2}) holds.
\end{definition}
To construct the desired sub- and supersolutions, we need some
preliminary lemmas.
\begin{lemma}\label{Le:2.1}
Let $V(x)$ satisfy $\rm(H_1)-(H_2)$, $\Lambda$ is defined by
(H$_3$). If $\Lambda<V_\infty$, then $\Lambda$ has a minimizer
$\varphi(x)\in C^{2,\gamma}_{loc}\cap W^{2,q}(\mathbb{R}^3)$ for
any $q\in(1,+\infty)$ with
\begin{equation}\label{eq:2.4}
-\Delta\varphi(x)+V(x)\varphi(x)=\Lambda\varphi(x), \ \
x\in\mathbb{R}^3,
\end{equation}
\begin{equation}\label{eq:2.5}
0<\varphi(x)<C|\varphi|_{\infty}e^{-l|x|},\ \ x\in\mathbb{R}^3,
\end{equation}
where $l\in(0,\sqrt{V_\infty-\Lambda})$ and $C=C(l,\Lambda)>0$ is a
constant.
\end{lemma}
{\it Proof.} It follows from Theorems 3.19 and 3.20 of
\cite{C.A.Stuart-NFaade} that there exists  $\varphi(x)\in
C(\mathbb{R}^3)\cap W^{2,2}(\mathbb{R}^3)$ such that
$(\ref{eq:2.4})$ $(\ref{eq:2.5})$ hold and $\varphi\in
L^1(\mathbb{R}^3)\cap L^\infty(\mathbb{R}^3)$ by (\ref{eq:2.5}).
Hence $\rm(H_1)$ and our Lemma \ref{Le:3.1} in Section 3 show that
$\varphi(x)\in W^{2,q}(\mathbb{R}^3)$ for any $q\in(1,+\infty)$.
This and embedding theorem implies that $\varphi(x)\in
C^{1,\gamma}_{loc}(\mathbb{R}^3)$. Thus, Theorem 9.19 in
\cite{DGilbargNStrudinger} gives that $\varphi(x)\in
C^{2,\gamma}_{loc}(\mathbb{R}^3)$.

\begin{lemma}\label{Le:2.2}
For any measurable function $u(x)$ on $\mathbb{R}^3$ with
 $$0<
u(x)\leqslant \psi(x)=\frac{1}{(1+|x|^2)^\alpha},  \mbox{ with }
\alpha>\frac{3}{4}.$$
  Then $\psi(x)\in L^q(\mathbb{R}^3)$ for all
$q\in[2,+\infty]$. Let
$\phi(x)=\int_{\mathbb{R}^{3}}\frac{u^{2}(y)}{4\pi|x-y|}dy$, we
have that
\begin{equation}\label{eq:2.6}
0<\phi(x)\leqslant\frac{1}{2(2\alpha-1)},\ \text{for any}\
x\in\mathbb{R}^{3}.
\end{equation}
\end{lemma}
{\it Proof.} By the definition, $\psi(x)\in
L^{\infty}(\mathbb{R}^3)$. Since $\alpha>\frac{3}{4}$ and
$2q\alpha
>3$ for $q\in[2, +\infty)$, it is not difficult to see that $\psi(x)\in
L^q(\mathbb{R}^3)$ for $q\in[2,+\infty]$.

 Let
$g(x):=\frac{1}{4\pi}\int_{\mathbb{R}^3}\frac{\psi^{2}(y)}{|x-y|}dy$.
Theorem 9.9 of \cite{DGilbargNStrudinger} shows that $g\in
D^{2,q}(\mathbb{R}^3)$  and hence $g(x)\in C^2(\mathbb{R}^3)$ by
Theorem 9.19 of \cite{DGilbargNStrudinger}, and
\begin{equation}\label{eq:2.6.1}
-\Delta g(x)=\psi^2(x),\ \ x\in\mathbb{R}^3.
\end{equation}
Since $\psi(x)$ is radially symmetric, it follows from
\cite[Proposition 4]{EdelsonALStuartCA-JDE} that $g(x)$ is
radically symmetric. Let $g(r)=g(|x|)$ and $\psi(r)=\psi(|x|)$ and
(\ref{eq:2.6.1}) becomes
\begin{equation}\nonumber
(r^2g'(r))'=-r^2\psi^2(r),
\end{equation}
then integrating over [0,r], we see that
\begin{equation}\nonumber
r^2g'(r)=-\int_0^{r}s^2\psi^2(s)ds<0.
\end{equation}
 This shows that $g(|x|)$ is strict decreasing. Hence, for any $x\in \mathbb{R}^3$
\begin{eqnarray}
\phi(x)&\leqslant&\frac{1}{4\pi}\int_{\mathbb{R}^{3}}\frac{\psi^{2}(y)}{|x-y|}dy
\leqslant
g(0)=\frac{1}{4\pi}\int_{\mathbb{R}^{3}}\frac{\psi^{2}(y)}{|y|}dy
 = \frac{1}{4\pi}\int_{\mathbb{R}^{3}}\frac{1}{|y|(1+|y|^2)^{2\alpha}}dy\nonumber\\
&=&\int_{\mathbb{R}^{3}}\frac{r}{(1+r^2)^{2\alpha}}dr
 =\frac{1}{2(2\alpha-1)}.\nonumber
\end{eqnarray}
So, (\ref{eq:2.6}) is proved. $\Box$

The following lemma gives a pair of sub- and supersolutions of
(\ref{eq:1.1}).
\begin{lemma}\label{Le:2.3} Under the assumptions of Theorem
\ref{th:1.1}, let $\varphi(x)$ and $\psi(x)$ be given by Lemmas
\ref{Le:2.1} and \ref{Le:2.2}, respectively. Then, for each
$\lambda\in[0,-2(2\alpha-1)\Lambda)$, there exists
$\epsilon_0\in(0,1)$ such that $\psi(x)$ and
$\epsilon_0\varphi(x)$ are  super- and  subsolutions of
(\ref{eq:1.1}), respectively. Moreover
$\epsilon_0\varphi(x)<\psi(x)$ for any $x\in \mathbb{R}^3$.
\end{lemma}
{\it Proof:}\quad For any $u\in W^{2,2}(\mathbb{R}^3)$, let
$\phi(x)=\int_{\mathbb{R}^{3}}\frac{u^{2}(y)}{4\pi|x-y|}dy$. Since
$u\in W^{2,2}(\mathbb{R}^3)\subset L^{q}(\mathbb{R}^3)\cap
C^{0,\gamma_1}(\mathbb{R}^3)$ for any $q\in[2,+\infty]$ and some
$\gamma_1\in(0,1)$, and the Hardy-Littlewood-Sobolev inequality
(see e.g. \cite[Theorem 4.3]{EHLiebMLoss-analysis}) yields that
$|\phi|_6\leqslant C|u|_{\frac{12}{5}}^2$ and $\phi\in
L^1_{loc}(\mathbb{R}^3)$. Then, it follows from \cite[Theorem
9.9]{DGilbargNStrudinger} that
\[
|D^2\phi|_q\leqslant C|u|_{2q}^2 \ \text{for each
}q\in(1,+\infty), \mbox{ and }
 -\Delta\phi(x)=u^2(x)\ a.e.\
x\in\mathbb{R}^3.
\]
So $u\in C^{0,\gamma_1}(\mathbb{R}^3)\cap L^q (\mathbb{R}^3)$ by
Sobolev embedding. Hence, Theorems 9.19 and 9.20 of
\cite{DGilbargNStrudinger} imply that $\phi\in
C_{loc}^{2,\gamma_1}(\mathbb{R}^3)$ and
\begin{equation}\label{eq:2.7.0}
-\Delta\phi(x)=u^2(x), \ \ x\in\mathbb{R}^3,\
\lim\limits_{|x|\rightarrow+\infty}\phi(x)=0.
\end{equation}
On the other hand, the uniqueness of the solution of
(\ref{eq:2.7.0}) implies that any solution of (\ref{eq:2.7.0})
with $u\in W^{2,2}(\mathbb{R}^3)$ must have the form of
$$\phi(x)=\int_{\mathbb{R}^{3}}\frac{u^{2}(y)}{4\pi|x-y|}dy.$$
Since $\rm(H_4)$ and $\Delta
\psi(x)=\frac{2\alpha[(2\alpha-1)r^2-3]}{(r^2+1)^2}\psi(x)=H(r)\psi(x)$,
it follows that
\begin{eqnarray}\label{eq:2.7}
\Delta \psi+Q(x){\psi}^{p}=H(r)\psi+Q(x){\psi}^{p} \leq
H(r)\psi+(V(x)-H(r)) \psi = V(x)\psi.
\end{eqnarray}
This and $\phi(x) \geq 0$ follow that, for each
$\lambda\geqslant0$,
\begin{equation}\nonumber
-\Delta\psi(x)+V(x)\psi(x)+\lambda\phi(x)\psi(x)\geqslant
Q(x)\psi^{p}(x),\ \text{on}\ \mathbb{R}^3,
\end{equation}
where $\phi(x)$ satisfies (\ref{eq:2.7.0}). So, $\psi(x)$ is a
supersolution of (\ref{eq:1.1}).

 For any $\epsilon>0$, since $\varphi$ satisfies (\ref{eq:2.4}), this
 yields
\begin{equation}\label{eq:2.9}
-\Delta\epsilon\varphi(x)+V(x)\epsilon\varphi(x)=\Lambda\epsilon\varphi(x)\
\text{for all}\ x\in\mathbb{R}^3.
\end{equation}
For any $\phi$ satisfying (\ref{eq:2.7.0}) with
$0<u(x)\leqslant\psi(x)$, if $\lambda\in[0,-2(2\alpha-1)\Lambda)$,
it follows from (\ref{eq:2.6}) that
\begin{equation}\label{eq:2.11}
\lambda\phi(x)+\Lambda\leqslant\frac{1}{2(2\alpha-1)}[\lambda+2(2\alpha-1)\Lambda]:=\delta_\lambda<0,\
\text{for any}\ x\in\mathbb{R}^{3}.
\end{equation}
By $\rm(H_1)$ and (\ref{eq:2.5}), $Q(x)\varphi^{p-1}\in
L^{\infty}(\mathbb{R}^3)$, then there exists a constant
$M\in(0,+\infty)$ such that
\begin{equation}\label{eq:2.12}
Q(x)\varphi^{p-1}\geqslant-M\ \text{ for any}\ x\in\mathbb{R}^{3}.
\end{equation}
From (\ref{eq:2.11}) and (\ref{eq:2.12}), there exists
$\epsilon_\lambda>0$ such that, for any
$\epsilon\in(0,\epsilon_\lambda)$
\begin{equation}\label{eq:2.13}
\lambda\phi(x)+\Lambda\leqslant\delta_\lambda
\leqslant-\epsilon^{p-1}M\leqslant\epsilon^{p-1}Q(x)\varphi^{p-1},\
\text{for all}\ x\in\mathbb{R}^{3}.
\end{equation}
 Then for each
$\lambda\in[0,-2(2\alpha-1)\Lambda)$, and
$\epsilon\in(0,\epsilon_\lambda)$, it follows from (\ref{eq:2.9}),
(\ref{eq:2.13}) and $\varphi(x)>0$ that
\[
-\Delta\epsilon\varphi(x) +
V(x)\epsilon\varphi(x)+\lambda\phi(x)\epsilon\varphi(x) =
(\Lambda+\lambda\phi(x))\epsilon\varphi(x) \leqslant
Q(x)(\epsilon\psi)^{p}(x), \  \text{ on }\ \mathbb{R}^3.
\]
 This means that $\epsilon\varphi(x)$ is a subsolution if $\epsilon\in(0,\epsilon_\lambda)$.
 Moreover, by (\ref{eq:2.5}) and the definition of
 $\psi$, we know that there exits $\epsilon_0\in(0,\epsilon_\lambda)$ such that
$\epsilon_0\varphi(x)<\psi(x)$ for any $x\in \mathbb{R}^3$.
 $\Box$

\section{Proofs of the main Theorems}
Now, we turn to showing our main Theorems \ref{th:1.1} and
\ref{th:1.2}. To prove Theorem \ref{eq:1.1}, an iteration sequence
is required, and it can be obtained by the sub- and supersolutions
given by Lemma 2.3 as well as the following lemmas.
\begin{lemma} \cite[Proposition 1]{EdelsonALStuartCA-JDE}\label{Le:3.1} Consider
$k>0$,
\begin{description}
 \item $(\rm i)$  For each $f\in L^q(\mathbb{R}^3)$ with $q\in[1,+\infty]$,
 there exists an unique $u:=Tf\in L^{q}(\mathbb{R}^3)$ satisfying $-\Delta
 u+ku=f$ on $\mathbb{R}^3$ in the sense of distributions.
 \item $(\rm ii)$ Let $f\in L^q(\mathbb{R}^3)$ with
 $q\in(1,+\infty)$, then $Tf\in W^{2,q}(\mathbb{R}^3)$. Moreover,
 there exists a constant $C=C(k,q)>0$ such that
$ \|Tf\|_{2,q}\leqslant C|f|_q. \Box $
\end{description}
\end{lemma}
\begin{lemma}\label{Le:3.2} Let $\epsilon_0\varphi(x)$ and $\psi(x)$
be given by Lemma \ref{Le:2.3}. Consider the following problem
\begin{equation}\label{eq:3.1}
-\Delta u(x)+ku(x)=f(x,w,v),\ \ x\in\mathbb{R}^3,
\end{equation}
where $k$ is a positive constant, $w$ and $v$ are functions on
$\mathbb{R}^3$, $f$:
$\mathbb{R}^3\times\mathbb{R}\times\mathbb{R}\rightarrow
\mathbb{R}$. For $\gamma\in(0,1)$ and
 $q\in(1,+\infty)$, we assume that
\begin{description}
 \item $(\rm F_1)$  $f(x,w,v)\in C^{0,\gamma}_{loc}(\mathbb{R}^3)\cap
 L^{q}(\mathbb{R}^3)$ if $w,v\in C^{2,\gamma}_{loc}(\mathbb{R}^3)\cap
 W^{2,q}(\mathbb{R}^3)$ and $\epsilon_0\varphi(x)\leqslant w,v\leqslant \psi(x)$.
 \item $(\rm F_2)$ $ f(x,u(x),\epsilon_0\varphi(x))\leqslant
f(x,u(x),u(x))\leqslant f(x,u(x),\psi(x))$ for any $u(x)\in
C^{2,\gamma}_{loc}(\mathbb{R}^3)\cap
 W^{2,q}(\mathbb{R}^3)$ with
$\epsilon_0\varphi(x)\leqslant u(x)\leqslant \psi(x).$
 \item $\rm ( F_3)$
$ -\Delta \epsilon_0\varphi(x)+k\epsilon_0\varphi(x)\leqslant
f(x,u(x),\epsilon_0\varphi(x))$ and $ -\Delta
\psi(x)+k\psi(x)\geqslant f(x,u(x),\psi(x)) $ for any $u(x)\in
C^{2,\gamma}_{loc}(\mathbb{R}^3)\cap
 W^{2,q}(\mathbb{R}^3)$ with
$\epsilon_0\varphi(x)\leqslant u(x)\leqslant \psi(x).$
\end{description}
Then, there exists $\{u_n\}\subset
C^{2,\gamma}_{loc}(\mathbb{R}^3)\cap
 W^{2,q}(\mathbb{R}^3)$ for any $q\in(1,+\infty)$ such that
 \begin{equation}\label{eq:3.2}
-\Delta u_{n+1}(x)+ku_{n+1}(x)=f(x,u_{n},u_{n}), \quad
x\in{\mathbb{R}}^{3},
\end{equation}
 \begin{equation}\label{eq:3.3}
\epsilon_0\varphi(x)\leqslant u_n(x)\leqslant \psi(x), \quad
x\in{\mathbb{R}}^{3},
\end{equation}
 \begin{equation}\label{eq:3.4}
\|u_{n+1}\|_{2,q}\leqslant C(k,q) |f(x,u_n,u_n)|_{q}.
\end{equation}
\end{lemma}
{\it Proof:}\quad  Let $u_{0}=\epsilon_0\varphi$ and Lemma
\ref{Le:2.1} implies that $u_0\in
C^{2,\gamma}_{loc}(\mathbb{R}^3)\cap
 W^{2,q}(\mathbb{R}^3)$ for any $q\in(1,+\infty)$. Then $f(x,u_0,u_0)\in
C^{0,\gamma}_{loc}(\mathbb{R}^3)\cap L^{q}(\mathbb{R}^3)$ by
$\rm(F_1)$. Applying Lemma \ref{Le:3.1} to problem (\ref{eq:3.1})
with $w=v=u_0$, we get $u_{1}(x)\in W^{2,q}(\mathbb{R}^3)$ such
that
\begin{equation}\label{eq:3.5}
-\Delta u_{1}(x)+ku_{1}(x)=f(x,u_{0},u_{0}), \quad
x\in{\mathbb{R}}^{3},
\end{equation}
 and then $u_{1}\in
C^{2,\gamma}_{loc}(\mathbb{R}^3)$ by Theorem 9.19 in
\cite{DGilbargNStrudinger}. Taking $u=u_0$ in $\rm(F_2)$ $(F_3)$
and noting that $\epsilon_0\varphi<\psi$, we see that
$$ f(x,u_0(x),\epsilon_0\varphi(x))\leqslant
f(x,u_0(x),u_0(x))\leqslant f(x,u_0(x),\psi(x)),$$
\[
-\Delta \epsilon_0\varphi(x)+k\epsilon_0\varphi(x)\leqslant
f(x,u_0(x),\epsilon_0\varphi(x)), \quad\quad x\in{\mathbb{R}}^{3},
\]
\[
-\Delta \psi(x)+k\psi(x)\geqslant f(x,u_0(x),\psi(x)), \quad\quad
x\in{\mathbb{R}}^{3}.
\]
These and (\ref{eq:3.5}) give that
\begin{equation} \label{eq:3.50}
-\Delta(\epsilon_0\varphi-u_{1})+k(\epsilon_0\varphi-u_{1})\leqslant0\leqslant
-\Delta(\psi-u_{1})+k(\psi-u_{1}) .
\end{equation}
Hence the maximum principle implies that
$
\epsilon_0\varphi(x)\leqslant u_{1}(x)\leqslant\psi(x).
$
On the other hand, Lemma \ref{Le:3.1} $\rm(ii)$ shows that
\begin{equation}\label{eq:3.51}
\|u_{1}\|_{2,q}\leqslant C(k,q) |f(x,u_0,u_0)|_{q}.
\end{equation}
Inductively, given $u_n\in C^{2,\gamma}_{loc}(\mathbb{R}^3)\cap
 W^{2,q}(\mathbb{R}^3)$ ($n=1,2,\cdots$) with $\epsilon_0\varphi(x)\leqslant u_n(x)\leqslant \psi(x)$,
 by Lemma \ref{Le:3.1} $\rm(i)$ and Theorem 9.19
of \cite{DGilbargNStrudinger}, there exists $u_{n+1}\in
C^{2,\gamma}_{loc}(\mathbb{R}^3)\cap
 W^{2,q}(\mathbb{R}^3)$  such that
\begin{equation}\label{eq:3.6}
-\Delta u_{n+1}(x)+ku_{n+1}(x)=f(x,u_{n},u_{n}), \quad
x\in{\mathbb{R}}^{3}.
\end{equation}
Taking $u=u_n$ in $\rm(F_2)$ and $\rm(F_3)$, similar to the
discussion of (\ref{eq:3.50}) and (\ref{eq:3.51}),  it follows
from (\ref{eq:3.6}) and the maximum principle  that $
\epsilon_0\varphi(x)\leqslant u_{n+1}(x)\leqslant\psi(x). $ Then,
by Lemma \ref{Le:3.1} $\rm(ii)$, $ \|u_{n+1}\|_{2,q}\leqslant
C(k,q) |f(x,u_n,u_n)|_{q}. $
 $\Box$

{\bf Proof of Theorem 1.1:} \quad For $v,w\in
W^{2,2}(\mathbb{R}^3)$, we denote
$$\phi_{w}(x)=\int_{\mathbb{R}^{3}}\frac{1}{4\pi|x-y|}w^{2}(y)dy,$$
 and for $\lambda\in[0,-2(2\alpha-1)\Lambda)$ and $k>0$ large enough, define
\begin{equation}\label{eq:3.7}
f(x,w,v)=Q(x)|v|^{p-1}u+kv-  V(x)v-\lambda\phi_{w}(x)v.
\end{equation}
 We
prove now the theorem by the following steps. In what follows,
$\epsilon_0\varphi$ and $\psi$ are the sub- and supersolutions given
by Lemma \ref{Le:3.2}.\\
{\bf Step 1: There exists $\{u_n\}\subset
C^{2,\gamma}_{loc}(\mathbb{R}^3)\cap
 W^{2,q}(\mathbb{R}^3)$ for any $q\in(1,+\infty)$
 such that

  \begin{equation}\label{eq:3.7.1}
-\Delta u_{n+1}(x)+ku_{n+1}(x)=f(x,u_{n},u_{n}), \quad
x\in{\mathbb{R}}^{3},
\end{equation}
 \begin{equation}\label{eq:3.7.2}
\epsilon_0\varphi(x)\leqslant u_n(x)\leqslant \psi(x), \quad
x\in{\mathbb{R}}^{3},
\end{equation}
 \begin{equation}\label{eq:3.7.3}
\|u_{n+1}\|_{2,q}\leqslant C(k,q) |f(x,u_n,u_n)|_{q},
\end{equation}
where $f$ is defined by (\ref{eq:3.7}).} \\
By Lemma \ref{Le:3.2}, Step 1 is proved if the function $f$
defined by (\ref{eq:3.7}) satisfies  $\rm(F_1)$ to $\rm(F_3)$ . By
 $\rm(H_1)$, Lemmas \ref{Le:2.2} and \ref{Le:2.3}, it is not
difficult to know that $\rm(F_1)$ and $\rm(F_3)$ hold. For
$0<w,v\leqslant\psi$ and $k>0$ large enough, it follows from $\rm
(H_1)$ and Lemma \ref{Le:2.2} that
\begin{equation}\label{eq:3.9}
\frac{\partial f(x,w,v)}{\partial v}
=pQ(x)|v|^{p-2}v+k-V(x)-\lambda\phi_{w}(x)>0,\mbox{ for any }
x\in\mathbb{R}^{3}.
\end{equation}
This implies that $\rm(F_2)$ holds. Hence, Step 1 is
complete.\\
{\bf Step 2: There exists $u\in W^{2,q}(\mathbb{R}^3)$ such that, by
passing to a subsequence, $\{u_n\}$ converges to $u$ weakly in
$W^{2,q}(\mathbb{R}^3)$ and strongly in $L^q(\mathbb{R}^3)$ for all
$q\in[2,+\infty)$.}\\
By Lemma \ref{Le:2.2}, $\psi\in L^2(\mathbb{R}^3)\cap
L^{\infty}(\mathbb{R}^3)$, then it follows from (\ref{eq:3.7}) and
$\rm(H_1)$ that
\begin{equation}\nonumber
|f(x,u_{n},u_{n})|_{{q}}\leqslant C|\psi|_{{q}},\ \text{for all }\
n\in\mathbb{N}\ \text{and }\ q\in[2,+\infty),
\end{equation}
where $C>0$ is a constant independent of $n$ and $\lambda$. So,
Lemma \ref{Le:3.2} implies that
\begin{equation}\label{eq:3.10}
\|u_{n}\|_{{2,q}}\leqslant C(k,q)|f(x,u_n,u_n)|_{{q}}\leqslant
C|\psi|_{{q}},\  q\in[2,+\infty).
\end{equation}
This means that $\{u_n\}$ is bounded in $W^{2,q}(\mathbb{R}^3)$ for
each $q\in[2,+\infty)$. So, passing to a subsequence, there is $u\in
W^{2,q}(\mathbb{R}^3)$ such that
\begin{equation}\nonumber
u_{n}\overset{n}{\rightharpoonup}u \ \text{weakly in }
W^{2,q}(\mathbb{R}^3)\ \text{ and }u_{n}\overset{n}{\rightarrow}u,\
a.e \text{ on } x\in\mathbb{R}^3.
\end{equation}
Therefore, $u\in\cap_{2\leqslant q<+\infty}W^{2,q}(\mathbb{R}^3)$.
By Lemma \ref{Le:2.2},
$\psi\in\cap_{q=2}^{+\infty}L^q(\mathbb{R}^3)$, then
$0<u_n^q(x)\leqslant\psi^q(x)$ by (\ref{eq:3.7.2}). Thus, the
dominated convergence theorem shows that
$$\int_{\mathbb{R}^3}|u_n(x)|^qdx\overset{n}{\rightarrow}\int_{\mathbb{R}^3}
|u(x)|^qdx,\ \text{for all }q\in[2,+\infty),$$ and Lemma 1.32 in
\cite{MWillem-Minitheorems} implies that
$|u_n(x)-u(x)|_q\overset{n}{\rightarrow}0.$

{\bf Step 3: $u\in C^2(\mathbb{R}^3)$ is a solution of (\ref{eq:1.1}).} \\
Multiplying (\ref{eq:3.7.1}) by $\eta(x)\in
C^{\infty}_0(\mathbb{R}^3)$, then integrating by parts over
$\mathbb{R}^3$, it yields that
\begin{eqnarray}\label{eq:3.12}
\int_{\mathbb{R}^{3}}\{-\Delta
u_{n+1}(x)+ku_{n+1}(x)\}\eta(x)dx&=&\int_{\mathbb{R}^{3}}u_{n+1}(x)\{-\Delta
\eta(x)+k\eta(x)\}dx\nonumber\\
&=&\int_{\mathbb{R}^{3}}f(x,u_{n},u_{n})\eta(x)dx,
\end{eqnarray}
 Since $u_n\overset{n}{\rightarrow}u\ a.e.\ \text{on}\
\mathbb{R}^3$, noting (\ref{eq:3.7.2}) and $\rm(H_1)$,  the
dominated convergence theorem shows that
\begin{eqnarray}\label{eq:3.13}
\lim_{n\rightarrow\infty}\int
_{\mathbb{R}^{3}}u_{n}(x)[-\Delta\eta(x)+k \eta(x)]dx =\int
_{\mathbb{R}^{3}}u(x)[-\Delta \eta(x)+k \eta(x)]dx,
\end{eqnarray}
\begin{eqnarray}\label{eq:3.14}
\lim_{n\rightarrow\infty}\int
_{\mathbb{R}^{3}}u_{n}(x)[k-V(x)-Q(x)u_{n}^{p-1}]\eta(x)dx =\int
_{\mathbb{R}^{3}}u(x)[k-V(x)-Q(x)u^{p-1}]\eta(x)dx.
\end{eqnarray}
 Letting $\phi_{n}=\phi_{u_n}$.
By Step 2, $|u_n(x)-u(x)|_q\overset{n}{\rightarrow}0$ for $q \in
[2,+\infty)$, it follows from \cite[Lemma 2.1]{RuizD-JFA}  that
 $\phi_n \stackrel{n} \rightarrow \phi_u $ in
$D^{1,2}(\mathbb{R}^3)$ and hence in $L^{6}(\mathbb{R}^3)$. For
any $\eta(x)\in C^{\infty}_{0}(\mathbb{R}^3)$, noting that $\|\eta
u_{n}\|_{L^{\frac{6}{5}}}$ is bounded, it follows from the
H\"older inequality and Sobolev embedding that
\begin{eqnarray}\label{eq:3.16}
&&\lim_{n\rightarrow\infty}\left|\int
_{\mathbb{R}^{3}}[\phi_{n}(x)u_{n}(x)\eta(x) -
\phi_{u}(x)u(x)\eta(x)]dx\right| \nonumber \\
& \leqslant & \lim_{n\rightarrow\infty}\left\{\int
_{\mathbb{R}^{3}}|\phi_{n}-\phi_{u}||\eta|u_{n}dx+\int
_{\mathbb{R}^{3}}\phi_{u}|\eta||u_{n}-u|dx\right\}\nonumber\\
&\leqslant&\lim_{n\rightarrow\infty}|\phi_{n}-\phi_{u}|_{{6}}|\eta
u_{n}|_{{\frac{6}{5}}}+o(1)
 \leqslant \lim_{n\rightarrow\infty}C|\nabla(\phi_{n}-\phi_{u})|_{2}|\eta
u_{n}|_{{\frac{6}{5}}}  .
\end{eqnarray}
By (\ref{eq:3.14}) (\ref{eq:3.16}) and the definition of $f$
(\ref{eq:3.7}), we see that
\begin{equation}\label{eq:3.17}
\lim_{n\rightarrow\infty}\int
_{\mathbb{R}^{3}}\{f(x,u_{n},u_{n})-f(x,u,u)\}\eta(x)dx=0, \mbox{
for any } \eta(x)\in C^{\infty}_{0}(\mathbb{R}^3).
\end{equation}
Thus, (\ref{eq:3.12}), (\ref{eq:3.13}) and (\ref{eq:3.17}) yield
 \begin{equation}\label{eq:3.17.1}
\int_{\mathbb{R}^{3}}\{-\Delta
u(x)+ku(x)\}\eta(x)dx=\int_{\mathbb{R}^{3}}u(x)\{-\Delta
\eta(x)+k\eta(x)\}dx=\int_{\mathbb{R}^{3}}f(x,u,u)\eta(x)dx,
\end{equation}
for any $\eta(x)\in C^{\infty}_{0}(\mathbb{R}^3)$. By the definition
of $f$, it gives that
\begin{equation}\label{eq:3.18}
\int_{\mathbb{R}^{3}}\{-\Delta
u(x)+V(x)u(x)\}\eta(x)dx+\int_{\mathbb{R}^{3}}\lambda\phi_u(x)u(x)\eta(x)dx
=\int_{\mathbb{R}^{3}}Q(x)u^{p}\eta(x)dx.
\end{equation}
 Since
$u\in\cap_{2\leqslant q<+\infty}W^{2,q}(\mathbb{R}^3)$, and the
embedding theorem shows that $u\in C^{1,\gamma}(\mathbb{R}^3)$. By
($\rm H_1$) and (\ref{eq:3.7}), $f(x,u,u)\in
C^{0,\gamma}_{loc}(\mathbb{R}^3)$. Then (\ref{eq:3.17.1}) and the
Theorem 9.19 of \cite{DGilbargNStrudinger} imply $u\in
C_{loc}^{2,\gamma}(\mathbb{R}^3)$.

 {\bf Step 4:
$0<u(x)\leqslant\psi(x)$ and $\|u\|_{2,q}\leqslant C$, where $C>0$
is a constant and independent of $\lambda $.}

 By (\ref{eq:3.10})
and the weakly lower semicontinuity of $\|\cdot\|_{2,q}$, we see
that
$
\|u\|_{2,q}\leqslant C|\psi|_{q},
$
where $C$ is a constant and independent of $\lambda $. By
$u_n\overset{n}{\rightarrow}u$ a.e. $\mathbb{R}^3$, and
$0<u_n(x)\leqslant\psi(x)$, we have $0<u(x)\leqslant\psi(x)$
a.e $\mathbb{R}^3$.$\Box$\\

{\bf Proof of Theorem 1.2:} By Theorem \ref{eq:1.1}, for each
$\lambda\in(0,-2(2\alpha-1)\Lambda)$, (\ref{eq:1.1}) has a solution
$u_\lambda$ satisfying (\ref{eq:3.18}), and
\begin{equation}\nonumber
\|u_{\lambda}\|_{2,q}\leqslant C, \ \text{for each }
q\in[2,+\infty),
\end{equation}
where $C$ is a constant independent of $\lambda$. Then, passing to a
subsequence, there exists $u_0\in W^{2,q}(\mathbb{R}^3)$ such that
\begin{equation}\nonumber
u_{\lambda}\overset{\lambda\rightarrow0}{\rightharpoonup}u_0, \
\text{weakly in }\ W^{2,q}(\mathbb{R}^3).
\end{equation}
Similar to Step 2 in the proof of Theorem 1.1, it follows from the
dominated convergence theorem and Lemma 1.32 of
\cite{MWillem-Minitheorems} that
\begin{equation}\label{eq:3.21}
\int_{\mathbb{R}^3}|u_{\lambda}(x)-u_0(x)|^qdx\overset{n}{\rightarrow}0.
\end{equation}
Finally, as Step 3 in the proof of Theorem 1.1, we know that, for
any $\eta(x)\in C^{\infty}_0(\mathbb{R}^3)$,
\begin{equation}\label{eq:3.22}
\int_{\mathbb{R}^{3}}\{-\Delta u_{0}(x)+V(x)u_{0}(x)\}\eta(x)dx
=\int_{\mathbb{R}^{3}}Q(x)u_0^{p}\eta(x)dx,
\end{equation}
and
 $u_0\in C_{loc}^{2,\gamma}(\mathbb{R}^3)\cap
W^{2,q}(\mathbb{R}^3)$ for all $q\in[2,+\infty)$, so $u_0$ is a
classical solution of (\ref{eq:1.2}). Moreover, using
(\ref{eq:3.18}) to (\ref{eq:3.22}), we know that
\begin{equation}\nonumber
\|u_{\lambda}-u_0\|_{W^{1,2}(\mathbb{R}^3)}\longrightarrow0,\
\text{as}\ \lambda\rightarrow0. \ \Box
\end{equation}


\noindent{\sc Wuhan Institute of Physics and Mathematics\\ Chinese
Academy of Sciences
\\P.O.Box 71710, Wuhan 430071, China\\}
Email:flymath@163.com and hszhou@wipm.ac.cn
\end{document}